\documentclass[a4paper,
               10pt,
               pdftex,
               normalheadings,
               headsepline,
               footsepline,
               headinclude,
               footinclude,
               DIV14,
               abstracton]
{scrartcl}

\usepackage[left=2.50cm, right=2.50cm, top=2.00cm, bottom=3.50cm]{geometry}

\usepackage
{
    graphicx,
    amssymb,
    amsmath,
    amsthm,
    xcolor,
    dsfont,
    algpseudocode,
    authblk,
}
\usepackage[ngerman, english]{babel}

\usepackage{enumerate}
\usepackage{algorithm}
\usepackage{tabularx}
\usepackage{subcaption} 
\usepackage{booktabs} 
\usepackage{siunitx}

\usepackage[colorlinks,
            pdffitwindow=false,
            plainpages=false,
            pdfpagelabels=true,
            pdfpagemode=UseOutlines,
            pdfpagelayout=SinglePage,
            bookmarks=false,
            colorlinks=true,
            hyperfootnotes=false,
            linkcolor=blue,
            urlcolor=blue!50!black,
            citecolor=green!50!black]
{hyperref}

\usepackage{sectsty}
\sectionfont{\fontsize{10}{10}\selectfont}
\subsectionfont{\fontsize{10}{10}\selectfont}

\DeclareMathOperator*{\argmax}{arg\,max}
\DeclareMathOperator*{\argmin}{arg\,min}

\newcommand*\diff{\mathop{}\!\mathrm{d}}
\renewcommand\vec{\mathbf}
\newcommand*\from{\colon}
\newcommand{\norm}[1]{\left\lVert#1\right\rVert}
\newcommand{\EE}[2]{\mathbb{E}_{#1\!\!}\left[#2\right]}
\newcommand{\CEE}[3]{\EE{#1}{{#2}~\middle\vert~{#3}}}

\newcommand{\tol}{\mathsf{tol}}

\let\OLDthebibliography\thebibliography
\renewcommand\thebibliography[1]{
	\OLDthebibliography{#1}
	\setlength{\parskip}{0pt}
	\setlength{\itemsep}{0pt plus 0.3ex}
}

\begin{document}
\title{Efficient time stepping for numerical integration using reinforcement learning}

\author[1]{Michael Dellnitz}
\author[2]{Eyke Hüllermeier}
\author[3]{Marvin Lücke}
\author[1]{Sina Ober-Blöbaum}
\author[1]{Christian Offen}
\author[4]{Sebastian Peitz}
\author[4]{Karlson Pfannschmidt}
\affil[1]{Department of Mathematics, Paderborn University}
\affil[2]{Department of Computer Science, LMU Munich}
\affil[3]{Modeling and Simulation of Complex Processes, Zuse Institute Berlin}
\affil[4]{Department of Computer Science, Paderborn University}

\date{}

\maketitle

\begin{abstract}
Many problems in science and engineering require an efficient numerical approximation of integrals or solutions to differential equations. For systems with rapidly changing dynamics, an equidistant discretization is often inadvisable as it either results in prohibitively large errors or computational effort. To this end, adaptive schemes, such as solvers based on Runge--Kutta pairs, have been developed which adapt the step size based on local error estimations at each step.
While the classical schemes apply very generally and are highly efficient on regular systems, they can behave sub-optimal when an inefficient step rejection mechanism is triggered by structurally complex systems such as chaotic systems.
To overcome these issues, we propose a method to tailor numerical schemes to the problem class at hand. This is achieved by combining simple, classical quadrature rules or ODE solvers with data-driven time-stepping controllers.
Compared with learning solution operators to ODEs directly, it generalises better to unseen initial data as our approach employs classical numerical schemes as base methods.
At the same time it can make use of identified structures of a problem class and, therefore, outperforms state-of-the-art adaptive schemes. Several examples demonstrate superior efficiency. Source code is available at \url{https://github.com/lueckem/quadrature-ML}.
\end{abstract}

\section{Introduction}
The numerical treatment of a vast number of problems in science and engineering requires the approximation of integrals, such as the computation of volume or mass flow, or the numerical solution of differential equations via \emph{Runge--Kutta methods} \cite{DB12}. Consequently, schemes for numerical discretization are a key element of scientific computing, and one of the main challenges is to determine a good trade-off between the required accuracy and numerical efficiency.

While the standard approach to developing \emph{quadrature rules} and \emph{numerical schemes to solve ODEs} is based on Taylor series expansions and the associated error bounds determined by higher-order derivatives \cite{DB12}, the advances in data science and machine learning have recently fueled the development of alternative concepts that are based on \emph{training data}. Most of these approaches are developed with the aim to efficiently compute the numerical solution of dynamical systems of high complexity, see, for instance, \cite{Piscopo2019,SIRIGNANO20181339,Liu2020HierarchicalDL,RPK17,RDQ19,KNP+20}, where the flow map $F$ that takes a state $x$ at time $t$ to a future state $x(t+\Delta t)$ is approximated.

In contrast to that, our work addresses the task of efficiently performing numerical integration for integrands or differential equations of a given problem class to a desired accuracy.
To this end, the next sample point at which the integrand or force term is evaluated is determined by finding an \emph{optimal trade-off} between the two conflicting criteria accuracy and numerical efficiency.
This task is carried out by a \emph{reinforcement learning} algorithm which, taking past function evaluations and learned knowledge about the problem class into account, determines the next sample point at which to evaluate the function of interest.
The efficiency of the proposed approach in comparison to state-of-the-art methods will be demonstrated using examples from the area of computing integrals (quadrature) as well as numerically solving differential equations.


The second key component of quadrature rules or numerical methods\,---\,the appropriate weighting of the different function evaluations\,---\,will be addressed briefly. 
While the classical Taylor-series-based construction aims to maximize the asymptotic order of a method on the class of analytic problems, given a finite accuracy and a restricted problem class, tailored weights can yield even more efficient schemes, as will be demonstrated by experiments.

We give a brief overview of numerical integration as well as reinforcement learning in Section~\ref{sec:preliminaries}.
In Section~\ref{sec:hierarchicalRL} we describe our method of learning an optimal time-stepping algorithm for specific problem classes from data, and how to optimize the weights of conventional quadrature rules in order to obtain an even more efficient scheme.
We apply these methods to multiple examples in Section~\ref{subsec:results} and compare the performance to state-of-the-art adaptive time-steppers.
Finally, we draw a conclusion in Section~\ref{sec:conclusion}.

\section{Preliminaries}\label{sec:preliminaries}

We will briefly review classical numerical quadrature rules and algorithms for differential equations as well as provide an overview of the reinforcement learning framework we will use in our algorithm.

\subsection{Numerical Quadrature}
The task in numerical quadrature is to efficiently determine the value of the integral 
\begin{equation} \label{eq:Integral}
	I=\int_{0}^{T} f(t) \mathrm{d}t
	\end{equation}
	from a few function evaluations, while maintaining a high accuracy.
	The most wide-spread strategy is to divide the interval $[0,T]$ into subintervals $I_i = [t_i,t_{i+1}]$ of length $h$ and to apply a quadrature formula 
	\begin{equation}\label{eq:quadrature}
	I_{i} = h \sum_{j=1}^s \omega_j f(t_i+c_j h)
\end{equation}
with fixed nodes $0 \le c_1<\ldots<c_s \le 1$ and weights $\omega_j \in \mathbb{R}$ in each subinterval \cite{Deuflhard2003}.
The standard numerical quadrature schemes implemented in popular software including MATLAB, Python's SciPy, or the FORTRAN subroutine package QUADPACK combine subdivision strategies with classical quadrature formulas. Their general purpose algorithms are most prominently based on Gauss-Kronrod quadrature \cite{scipyIntegrate,Quadpack,SHAMPINE2008131}.
Gauss-Kronrod quadrature approximates the value of the integral using $s=2r+1$ nodes and weights such that polynomials up to including order $3r+1$ can be integrated exactly (degree of exactness). The error is of order $3r+3$ in the length of the biggest subinterval $h$. A subset of $r$ nodes can be used in a Gauss quadrature formula such that another approximation with degree of exactness $2r-1$ and of order $2r+2$ in $h$ can be obtained without any further evaluations of the integrand. Comparing both approximations provides an estimate for the numerical error in the subinterval, based on which further subdivisions are performed.
The nodes are in an irrational relation which enhances the robustness of the error estimate, as argued in \cite{SHAMPINE2008131}, for instance. However, as a consequence, previous function evaluations cannot be used in subsequent steps.

Many software packages employ the above mentioned adaptive Gauss-Kronrod subdivision algorithm with 21 nodes ($r=10$) as a default, which we will refer to as \texttt{GK21}.
Moreover, we denote the associated Kronrod quadrature rule with 21 nodes by \texttt{K21} and the Gauss quadrature rule with 21 nodes by \texttt{G21}.


\subsection{Numerical integration of differential equations}\label{subsec:integration}

Consider the initial value problem $\dot{x} = f(t,x)$ with initial condition $x(0)=x_0$ on the time interval $[0, T]$. 
Reformulating the problem in integral form
\begin{equation}\label{eq:SolutionODE}
    x(T) = x_0 + \int_{0}^{T} f(t, x(t)) \mathrm{d}t
\end{equation}
and applying numerical quadrature rules gives rise to integration methods for ordinary differential equations such as \emph{Runge--Kutta} methods \cite{DB12}. Time-stepping proceeds as follows:
for a given step size $h_0>0$, first the approximate value $x_1$ of the solution $x(t_1)$ at time $t_1=t_0+h_0$ is computed. Next, the process is repeated with $(t_1,h_1,x_1)$ replacing $(t_0,h_0,x_0)$ to compute $x_2$ for $h_1>0$. Iteration yields an approximation $x_0,x_1,x_2,\ldots$ of the solution $x$ at times $t_0,t_1,t_2,\ldots$, where $t_i = t_{i-1}+h_{i-1}$ is recursively defined.

While in the simplest case all step sizes $h_i$ are identical (constant time-stepping), more elaborate schemes choose step sizes dynamically based on error estimates that are available from previous steps. In this way, numerical integration can benefit from large time steps in regions where the solution does not change rapidly while satisfying accuracy requirements in regions which require small step sizes \cite{Hairer1993,Hairer1996}. 
Efficient time-stepping methods that are implemented in MATLAB or Python’s SciPy package \cite{scipy} include most prominently the Dormand--Prince pair \texttt{RK45}, consisting of a fifth order Runge--Kutta scheme together with a fourth order method that is used to estimate the local error of each step. Based on the error estimate, a step size controller then computes the step size for the next iteration \cite{DORMAND1980,Shampine1997}.
If the error estimate of a step is larger than some desired tolerance, the proposed step is rejected and the time-stepper tries again with a smaller step size, hence, wasting function evaluations.
Other methods include the 8th order scheme \texttt{DOP853} \cite{Hairer1993,SciPyODE}, or multi-step strategies for stiff equations \cite{Hairer1993} as well as automatic solver selection strategies \cite{WolframNDSolve}. 

The above-mentioned black box solvers are designed for an application to a broad range of problems. However, for problems with additional structure (for instance, Hamiltonian systems), specialized algorithms achieve significant advantages over off-the-shelf methods by making use of symmetries, conserved quantities, or symplecticity \cite{GeomIntegration}. On the other hand, a lack of regularity can cause order reductions or instabilities of classical integration schemes that are hard to overcome, though approaches are available in some cases \cite{Schratz2019}.

\subsection{Data-driven approaches to time-stepping}\label{subsec:RL}

While specialized numerical integration schemes have traditionally been developed by manually identifying structures first and then designing appropriate numerical schemes, in the following we will explore a data-driven approach to develop time-stepping strategies that adapt automatically to given problem classes and error tolerances.
In our approach, these strategies are learned by a neural network for specific problem classes, such as irregular integrands or dynamical systems that exhibit chaotic motions. Here, classical adaptive schemes can suffer from inefficiencies, for instance because proposed steps get rejected frequently when an error estimate becomes too large.
In comparison, our data-driven scheme can learn features specific to the problem class and use that knowledge to optimize its time-stepping strategy.

\subsubsection{Reinforcement learning} \label{subsec:reinforcement_learning}
In reinforcement learning, an agent acting in an environment $E$, typically modelled in the form of a Markov decision process $(\mathcal{S}, \mathcal{A}, p, r)$, is considered with $\mathcal{S}$ the set of possible states, $\mathcal{A}$ the set of possible actions, $p$ the so-called transition function, and $r$ the reward function. 
At each discrete time step $i$, the agent finds itself in a state $s_i \in \mathcal{S}$ of the environment. 
It decides on an action $a_i \in \mathcal{A}$ and receives an immediate reward $r_i = r(s_i, a_i)\in \mathbb{R}$ that depends on the current state and the action taken (and perhaps also on the successor state).
The behavior of the agent is captured by its \emph{policy} $\pi \from \mathcal{S} \to P(\mathcal{A})$ that prescribes actions given states. More specifically, the policy is not necessarily deterministic and defines a probability distribution over the actions available in a given state.
Likewise, the environment $E$ can be stochastic: the successor state $s_{i+1}$ resulting from action $a_i$ in the current state $s_i$, is determined by the transition function $p:\, \mathcal{S} \times \mathcal{A} \to  P(\mathcal{S})$. In the following, we will use the notation $r_i, s_{i+1} \sim E$ as a shorthand whenever $s_i$ and $a_i$ are clear from the context.

Given a sequence of states and actions $(s_i, a_i, s_{i+1}, a_{i+1}, \dots, s_T)$ starting at state $s_i$,
$R_i$ denotes the sum of discounted future rewards
$R_i = \sum_{j=i}^{T-1} \gamma^{(j-i)} r(s_j, a_j)$,
where $\gamma \in [0, 1]$ is a discount factor.
The actions are the result of sampling from the (stochastic) policy $\pi$,
which is not known in advance.
The goal is to learn a policy which, given an initial state
$s_0 \in \mathcal{S}$, maximizes the expected sum of discounted future rewards.

We consider the expected future rewards with respect to a specific action. The action-value function is defined as follows:
\begin{equation}
    Q^{\pi}(s_i, a_i) = \CEE{\pi}{R_i}{s_i, a_i}
\end{equation}
It assumes that the agent picks action $a_i$ in state $s_i$, while subsequently
picking actions according to policy $\pi$.
This definition can be further unrolled into a recursive definition:
\begin{equation}
   Q^{\pi}(s_i, a_i) =
    \EE{r_i,s_{i+1} \sim E}{
        r(s_i, a_i) + \gamma
        \EE{a_{i+1} \sim \pi}{Q^{\pi}(s_{i+1}, a_{i+1}}
    }
\end{equation}
This equation is known as the Bellman equation \cite{Lillicrap2016, sutton2018}.
Notice that the outer or inner expectation is not needed if the environment or policy function is deterministic.
The function
\begin{align}
    Q(s_i, a_i) =
    \EE{r_i,s_{i+1} \sim E}{
        r(s_i, a_i) + \gamma \max_{a_{i+1} \in \mathcal{A}} {Q(s_{i+1}, a_{i+1})}
    },
\end{align}
evaluates the action $a_t$ in state $s_t$
based on the premise that optimal actions are chosen in subsequent steps,
and is called \emph{Q-function}.
The optimal policy $\pi^*\from \mathcal{S} \to \mathcal{A}$ can be characterized as
\begin{equation}
    \pi^*(s) = \argmax_{a\in \mathcal{A}} Q(s, a).
\end{equation}

In this article, we use a reinforcement learning approach based on Q-learning \cite{watkins1992q}, where the $Q$-function is approximated by a neural network parametrized by~$\theta$, which we will denote by $Q_{\theta}$.
The parameter $\theta$ is initialized using a suitable initialization scheme
\cite{glorot2010}.
The training then proceeds over the course of several episodes.
In each episode the learner starts in a state $s_0 \in \mathcal{S}$ of the environment and explores the state space using the function $Q_{\theta}$ in conjunction with a probabilistic policy.
This results in a trajectory $\left((s_i, a_i, r_i) \right)_{i=0}^H$, where
$s_i$ is the state in step $i$, $a_i$ the chosen action, and $r_i$ is the resulting reward.
The horizon $H$ is the final time step of the trajectory and marks the end of one episode.
The training targets $Q'$ are defined as follows:
\begin{align} \label{eq:Q_training}
    Q'(s_i, a_i) = r_i + \gamma \max_{a\in\mathcal{A}} Q_\theta(s_{i+1}, a), 
\end{align}
where $\gamma \in [0, 1]$ is a discount factor for future rewards.
Given a loss function $L\colon \mathbb{R} \times \mathbb{R} \to \mathbb{R}$, the parameter vector
$\hat{\theta}$, to be used in the next episode, is calculated using empirical risk minimization:
\begin{equation} \label{eq:Q_training_update}
    \hat{\theta} \gets \argmin_{\theta \in \Theta}\frac{1}{H+1}\sum_{i=0}^H L\bigl(Q'(s_i, a_i), Q_{\theta}(s_i, a_i)\bigr) + R(\theta),
\end{equation}
where $R(\theta)$ is an additional regularization term.
In the following, we use the typical $L_2$ loss, i.\,e.,
$L(y, \hat{y}) = \norm{y - \hat{y}}_2$.
We continue to run episodes until the Q-function converges.
This approach and similar variants are often referred to as \textit{deep Q-learning} \cite{Mnih2015}.
While convergence to the exact Q-function is only rigorously proved for classical Q-learning in a finite state and action space setting \cite{KS99}, there are numerous publications showing the success of deep Q-learning in many applications,
see \cite{pmlr-v120-yang20a} for an overview.

\section{An RL algorithm for efficient adaptive integration} \label{sec:hierarchicalRL}

Our algorithm addresses two main challenges of numerical integration schemes:
\begin{enumerate}[(i)]
    \item the selection of step sizes in situations of rapidly changing behavior (e.\,g., bursts),
    \item qualitative changes in behavior, for instance in hybrid or switched systems, or systems which chaotic and regular regions.
\end{enumerate}
The approach is to train a model that can recommend optimal step sizes to use during integration and is superior to classical adaptive schemes.
We define an optimal step size as being as large as possible while a certain desired error tolerance is not exceeded in the integration step.
This allows us, in particular, to address issues such as the use of too large step sizes (``overshooting'') when the integrand changes rapidly, while being too conservative in other regions.
The proposed algorithm employs reinforcement learning to adapt to a particular class of problems.
The learning only has to be conducted once, and the obtained model can then be used to integrate similar functions much more efficiently.
We describe the training of the efficient time-stepper in Section~\ref{subsec:basisLearner}.

Moreover, given a specific class of problems, the weights of classical quadrature and integration methods may be suboptimal. We show how to calculate weights that are specifically tailored to the problem class in order to obtain an even more efficient scheme in Section~\ref{subsec:optim_weights}.

\subsection{Training the time-stepper via reinforcement learning}\label{subsec:basisLearner}

We first discuss training an efficient time-stepper for quadrature tasks in Section~\ref{subsec:method_quad}. Then we focus on the numerical integration of ordinary differential equations (ODEs) in Section~\ref{subsec:method_ode}.

\subsubsection{Quadrature tasks} \label{subsec:method_quad}
We assume that the functions we want to integrate are sampled from a set $X$ with probability measure $P$, i.e., a class of problems, and aim to train a neural network (NN) to perform the subdivision of an initial integration interval into subintervals.
The quadrature rule that is used to integrate the subintervals is fixed.
We restrict the algorithm to integrate ``from left to right'', that is, given the current subinterval $[t_i, t_{i+1}]$, the NN chooses the next step size (subinterval length) $h^+$ so that the next subinterval is given by $[t_{i+1}, t_{i+1} + h^+] = [t_{i+1}, t_{i+2}]$.
Provided with a sufficient amount of training data, the NN may be able to predict the future course of the function to some extent and, thus, choose a suitable step size.

Assuming the quadrature rule uses $s$ function evaluations $(f(t_i + c_j h))_{j=1,\dots,s}$ in the subinterval $[t_i, t_{i+1}]$, where the $c_j$ are the nodes of the quadrature rule (cf. \eqref{eq:quadrature}), we provide the subinterval length $h$ as well as the function evaluations as an input to the time-stepper.
It can then select the subsequent subinterval length $h^+$ from a finite set of options $\{h_1, \dots, h_n\}$.
 

The time-stepping NN is trained via Q-learning \cite{watkins1992q}. More specifically, the NN approximates the Q-function, where
Q receives the \textit{state} $s_i$ at step $i$ and an \textit{action} $a_i$ as inputs which, in our case, are
\begin{align*}
	s_i &= (h, f(t_i + c_1 h), \dots, f(t_i + c_s h)) \, , \\
	a_i &= h^+ \, .
\end{align*}
As described in equation \eqref{eq:Q_training}, the value of the Q-function can be defined implicitly as
\begin{align*}
	Q(s_i, a_i) = r_i + \gamma \max_{a_{i+1} \in \mathcal{A}} Q(s_{i+1}, a_{i+1}),
\end{align*}
with reward $r_i$ and discount factor $\gamma \in [0,1]$.
In preliminary experiments, a discount factor of $\gamma = 0$
resulting in a myopic Q-function yielded the same performance as
non-zero values.
Therefore, we simply set $\gamma = 0$ in all experiments.
The neural network receives $s_i$ as an input and returns all possible Q-values, i.\,e., the output is defined as
\begin{align*}
	(Q(s_i, h_1), \dots, Q(s_i, h_n)) \, .
\end{align*}
We then simply select the action with the highest Q-value, i.\,e., the highest expected reward, as our next step size.
The training of the neural network is conducted as described in Section~\ref{subsec:reinforcement_learning}, that is, for every episode a function $f$ is sampled from the problem class, integrated by the time-stepper, and the resulting rewards used to optimize the weights of the NN. To ensure sufficient exploration during the training phase, we do not always select the step size with the highest Q-value, but employ a probabilistic policy that sometimes selects one of the other step sizes randomly.

\subsubsection*{Reward function}
To learn a strategy that selects a step size that is as large as possible while an error tolerance $\tol$ is not exceeded, we have to choose an appropriate reward function.
In our case, the reward $r_i$ should depend on the integration error $\varepsilon$ in the integration step $t$ and on the selected step size $h^+$.
Here, the integration error $\varepsilon$ is defined as
\begin{align*}
    \varepsilon = | I - \hat{I} |,
\end{align*}
where $I = \int_{t_i}^{t_{i+1}} f(t) \diff t$ is the exact integral of the current step and $\hat{I}$ denotes the result of the quadrature formula.
If the exact integral $I$ is not known, we approximate it numerically with high accuracy (e.\,g., using a cumulative quadrature with step size multiple orders smaller than $h^+$).
A straightforward definition of the reward function is
\begin{align*}
	r_i &= \begin{cases}
		0, & \varepsilon > \tol\\ 
		h^+, & \varepsilon < \tol
	\end{cases}.
\end{align*}
\begin{figure}[t]
    \centering
    \includegraphics{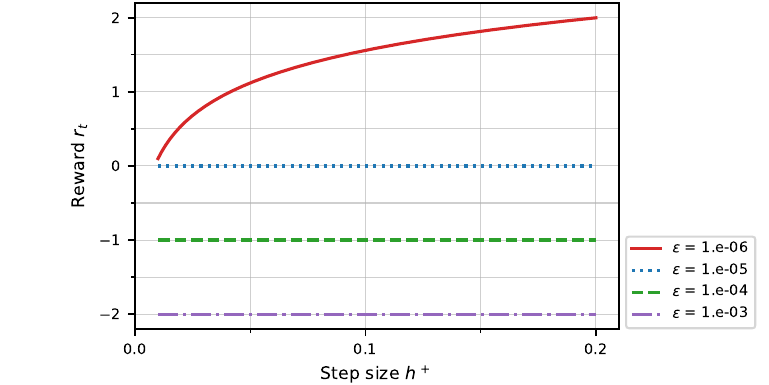}
    \caption{Behavior of the reward function \eqref{eq:rewardlog10} for varying $\epsilon$. The error tolerance is fixed at \num{1e-5}.}
    \label{fig:reward_function}
\end{figure}
However, in practice, it is advantageous to use a reward function with negative rewards for larger errors.
For problems exhibiting chaotic behavior such as the ones we consider in this paper, we expect integration errors to be scattered across several orders of magnitude.
This is why we chose to use the following reward function
\begin{align}\label{eq:rewardlog10}
	r_i &= \begin{cases}
	  \log_{10}(\frac{\tol}{\varepsilon}), & \varepsilon > \tol\\
	  a \log(b \cdot h^+), & \varepsilon < \tol
	\end{cases},
\end{align}
which combines both properties. It generates a reward of  $r_i = -m$ if $\varepsilon = 10^m \cdot \tol$.
The parameters $a$ and $b$ are chosen such that the positive and negative rewards are roughly on the same scale.
The behavior of reward function \eqref{eq:rewardlog10} with respect to
varying integration error $\epsilon$ is shown in Figure~\ref{fig:reward_function}.
This is the reward function we employ for all experiments throughout the paper.
Generally speaking, problem-depending reward functions can be necessary for efficient training.

\subsubsection*{Neural network architecture}
Throughout this paper, we use fully connected neural networks with four hidden layers and five times the number of input nodes per hidden layer. It should be noted that for the examples studied in the following, a significantly smaller network topology might be sufficient.
All NNs use rectified linear unit (ReLU) activation functions and employ the Adam stochastic gradient descent optimizer \cite{KB14} for training.

\subsubsection*{Limitations} It is apparent that the above constructed method can only be successful if the information contained in the state $s_i$, i.e., the function evaluations in one subinterval $[t_i, t_{i+1}]$, is sufficient to form some belief about the structure of $f(t)$ for $t>t_{i+1}$.
Thus, we do not expect the method to perform particularly well on very broad classes of problems, where past function evaluations are not informative enough about possible continuations of the function.
We do expect our method to excel on function classes of which the functions share some common and predictable behavior.

When the information contained in the state $s_i$ is not sufficient to form an adequate prediction of the structure of $f$, one can supplement past states to the input of the time-stepper, thereby increasing the available information. However, this considerably increases the complexity of the learning problem and, in the experiments that we conducted, did not lead to a significant increase in performance that would justify the additional complexity.

\subsubsection{Numerical Integration of Differential Equations} \label{subsec:method_ode}
The approach discussed above can readily be extended to the numerical integration of ordinary
differential equations (ODEs).
Here Runge-Kutta methods take the role of quadrature rules.
A Runge-Kutta method with $s$ stages can be written as
\begin{align*}
    x(t_{i+1})-x(t_i) = \int_{t_i}^{t_{i+1}} f(t, x(t)) \text{d}t \, \approx \, h \sum_{j=1}^s b_j k_j, \quad h = t_{i+1} - t_i \, .
\end{align*}
The stages $k_j$ correspond to approximations of the function values $f(\tilde{t}_j, x(\tilde{t}_j))$ for $s$ distinct time points $\tilde{t}_j \in [t_i, t_{i+1}]$, and thus, correspond to the function evaluations at the nodes in the quadrature setting.
The weights $b_j$ correspond to the weights of a quadrature rule.

Hence, we can apply the same method as presented in Section \ref{subsec:method_quad} in order to train an efficient time-stepper.
The input of the NN is now given by the state
\begin{align*}
    s_i = (h, k_1, \dots, k_s).
\end{align*}
Note that the stages $k_j$ are vectors, such that the total input dimension of the NN is $1 + s d$, where $d$ is the dimension of $x(t)$.
We define the integration error for the step $i$, which is needed to calculate the reward, as the deviation (in 2-norm) of the Runge Kutta estimate for $x(t_{i+1})$ from the correct solution when starting at the same initial value as the Runge Kutta estimate (local, step-wise error).
If the exact solution $x(t_{i+1})$ is not known, we approximate it numerically with sufficiently high accuracy.



\subsection{Optimization of weights}
\label{subsec:optim_weights}

Now that we have constructed a reinforcement learning approach, which is able to train accurate models for
efficient time-stepping, it is natural to ask what happens, if we adapt the weights of the quadrature rule to the given problem class as well.
Recall that a quadrature rule is of the following general form (cf. \eqref{eq:quadrature})
\begin{align*}
    \hat{I}_i = h \sum_{j=1}^s \omega_j f(t_i+c_j h)
\end{align*}
for suitable choices of the nodes $c_j$ and the weights $\omega_j$.
In case of ODEs, Runge-Kutta methods use the recurrence relation
\begin{align*}
    x_{i+1} = x_i + h \sum_{j=1}^{s} b_j k_j\ ,
\end{align*}
where the $b_j$ are the weights, $k_j$ the stages, and $x_i$ the states of the trajectory.

We propose to optimize these weights (i.\,e., $\omega_j$ or $b_j$) using linear regression, such that the expected integration error for the specific problem class is minimized.
For the case of numerical quadrature, we consider input instances of the form
$\vec{f}_i = (f(t_i+c_1 h), \dots, f(t_i+c_s h))$, while the ground-truth integral $I_i$ is the output to be predicted. (Again, if $I_i$ is not analytically known, we approximate it numerically with high accuracy.)
The optimization problem
\begin{align*}
    \vec{w}^* = \argmin_{\vec{w} \in \mathbb{R}^{s}} \sum_{i=1}^n ( h_i\vec{w}^T\vec{f}_i - I_i)^2,
\end{align*}
with $n$ being the size of the dataset, can then be solved using off-the-shelf libraries.

The application to Runge-Kutta methods is straightforward as well.
During each evaluation of our reinforcement learning model on a random trajectory, we record the vectors
$\vec{k} = (k_1, \dots, k_{\ell})$ to be used as input, while we
let $\Delta = x(t_{i+1}) - x_{i}$ be the corresponding output to be predicted, where $x(t_{i+1})$ is the correct solution when starting at $x(t_i) = x_i$.
(Again, if $x(t_{i+1})$ is not analytically known, we approximate it numerically with high accuracy.)
In case the state-space of the ODE is multi-dimensional (with $d$ denoting the number of dimensions), we split each instance into $d$ separate ones,
since we are only interested in one global set of weights.
As before, we solve
\begin{align*}
    \vec{b}^* = \argmin_{\vec{b} \in \mathbb{R}^{s}} \sum_{i=1}^n (h_i\vec{b}^T\vec{k}_i - \Delta_i)^2,
\end{align*}
to obtain an optimized weight vector $\vec{b}^*$.

As the optimal time-stepping technique may depend on the integration formula used, and vice-versa, it is recommended to train the time-stepping and to optimize the weights synchronously. To avoid the weights to be overly adapted to one particular trajectory,
we use an exponentially weighted moving average
\begin{align*}
    \vec{b}_{\mathrm{ep} + 1} = \alpha \vec{b}^* + (1 - \alpha) \vec{b}_{\mathrm{ep}}
\end{align*}
to update the weights, where $\mathrm{ep}$ was the last reinforcement learning epoch performed and $\vec{b}^*$ are the optimized weights for that epoch.
In our experiments $\alpha=0.05$ produced stable results.
In the following numerical experiments, we evaluate the performance of our reinforcement learning approach both with and without weight optimization.

\section{Numerical experiments}\label{subsec:results}

In this section, we will apply the methods presented above to several example problems and compare their performance to state-of-the-art adaptive techniques.
First, we will discuss learning an optimal subdivison of an interval into subintervals for challenging quadrature tasks in Section~\ref{subsec:ResultsQuadrature}.
Then we address learning an optimal time-stepper for the numerical solution of differential equations in Section~\ref{subsec:ResultsODE}. 

\subsection{Numerical quadrature}
\label{subsec:ResultsQuadrature}
The most challenging tasks in numerical quadrature are, besides functions with singularities or jumps, functions that exhibit a quickly changing and erratic behavior on different scales, so that the optimal subinterval sizes differ significantly.
We show that the efficient time-stepping method we propose can successfully deal with such functions in the following example.

\subsubsection{Example: Double Pendulum}

Due to their frequent and sharp changes of direction, it is difficult to perform quadrature on the angular velocities of a double pendulum, e.\,g., in order to recover the angular positions of the two pendulums.
The double pendulum system is explained in more detail in Section~\ref{subsec:example_doublependulum}, where we learn a time-stepper to efficiently solve the double pendulum ODE.
Here, we will assume a trajectory $f(t)$ of the angular velocity of the first pendulum is given for $0 \leq t \leq 100$ and the task is to calculate the integral $\int_0^{100} f(t) \text{d}t$.

We calculate the trajectory $f(t)$ using the \texttt{RK45} adaptive solver on the double pendulum ODE with very high accuracy. The initial condition for the double pendulum is randomly sampled from a fixed energy set.
Hence, we obtain a class of functions $f$ equipped with some probability measure, and the goal is to find an optimal efficient time-stepping technique to divide the integration interval into subintervals. Here, ``optimal'' refers to using the least expected function evaluations while maintaining an expected integration error lower than some bound (expectation with respect to the above mentioned probability measure).

We set the error tolerance to $\tol = 10^{-7}$ and train a time-stepping model (as described in Section~\ref{subsec:basisLearner}) with 20 step sizes evenly spaced on a log scale $h \in [0.1, 0.7]$.
For the quadrature itself we employ the Kronrod rule with 21 nodes (\texttt{K21}) on each subinterval for a fair comparison to the Gauss-Kronrod adaptive subdivision algorithm (\texttt{GK21}), which is the default adaptive quadrature rule in many software packages like Python's \textit{SciPy}.

\begin{figure}
    	\centering
    	\begin{subfigure}[b]{0.49\textwidth}
    	\includegraphics[width=\textwidth]{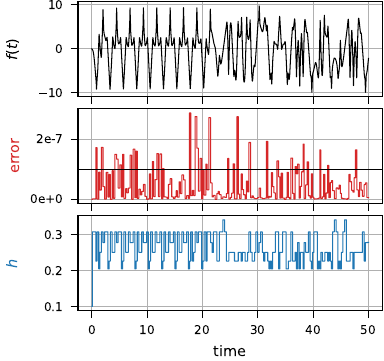}
    	\caption{time-stepping}
    	\label{fig:quad_a}
    	\end{subfigure}
    	\begin{subfigure}[b]{0.49\textwidth} \includegraphics[width=\textwidth]{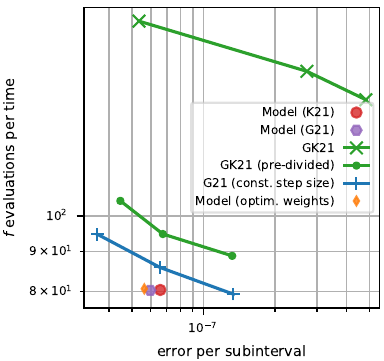}
    	\caption{Performance comparison}
    	\label{fig:quad_b}
    	\end{subfigure}
    	\caption{Quadrature of double pendulum. (a) time-stepping (division into subintervals) obtained by the model. (b) Performance comparison between the trained model and \texttt{GK21}. We show the model with \texttt{K21} and \texttt{G21} quadrature rules as well as optimized weights for the \texttt{G21} nodes. ``\texttt{GK21}'' refers to the adaptive subdivision algorithm applied to the whole integral from $t=0$ to $t=100$ (plotted are three different tolerances for subdivision); for ``\texttt{GK21} (pre-divided)'' the whole interval was first pre-divided into subintervals of length $0.35$ and then \texttt{GK21} was applied to each subinterval. Furthermore, we show \texttt{G21} applied to subintervals of constant sizes. All points are the average over an ensamble of 20 integrals, where the integrand $f(t)$ is sampled as described above.}
    	\label{fig:quad}
    \end{figure}

The results are shown in Figure~\ref{fig:quad}.
As the instance of a typical time-stepping in \ref{fig:quad_a} indicates, the model predominantly selects time steps (subinterval lenghts) between 0.2 and 0.3 and achieves subinterval integration errors close to the desired tolerance.
The average error, which is shown in \ref{fig:quad_b}, is slightly below the desired tolerance at about \num{7d-8} and the model uses approximately 80 function evaluations per time unit.
In comparison, the standard adaptive subdivision algorithm \texttt{GK21} needs substantially more evaluations to achieve a similar error. The reason is that it needs to execute many subdivisions before subintervals are of an appropriate length ($\sim 0.3$) for the required tolerance, and already evaluated function values are not reusable in the next subdivision step%
\footnote{For example, dividing the initial interval $[0,100]$ into $[0, 50]$ and $[50, 100]$ already ``costs'' 21 function evaluations, and so on.}.
Hence, this comparison is not fair, as we only gave interval lengths of appropriate size ($0.1 \leq h \leq 0.7$) to our model as options.
Therefore, we also computed the performance of \texttt{GK21} after having pre-divided the whole integration interval $[0,100]$ into subintervals of length 0.35, and then applying \texttt{GK21} to each subinterval.
As can be seen in Figure~\ref{fig:quad_b}, this approach still performs rather poorly compared to our model.

Note that there exist adaptive subdivision schemes, where function evaluations can be reused in subsequent subdivision steps, e.g., Romberg quadrature.
Such quadrature rules necessarily evaluate the integrand on nodes that are in a rational relation. However, high order quadrature rules, such as Gauss quadrature, require the evaluation of the integrand on an irrational grid. Moreover, it has been argued \cite{SHAMPINE2008131} that  using irrational grids improves the robustness of error estimations as correlations between function values are avoided.
\footnote{In interpreted languages such as Python or MATLAB the bottleneck to efficiency is in many cases the interpreter rather than evaluations of complex integrands. The default choices of quadrature algorithms have been justified in \cite{SHAMPINE2008131} by the fact that after each subdivision the positions, where the integrand needs to be evaluated, can be passed to the function representing the integrand all at once in a vectorized form. In this context, effectively only the number of times a subdivision needs to be performed counts.
}
Independently on whether function evaluations can be reused or not, subdivision techniques are generally not efficient for the task at hand.
The reason is that small tweaks to the step size are required, e.g., between $0.2$ and $0.3$, whereas subdivision techniques at least double the number of evaluations in a subdivision step.
Hence they often produce an unnecessarily precise approximation when they choose to further subdivide a subinterval, or a too large error when they choose to not subdivide.

Due to the inefficiency of subdivision techniques, we also compared our trained model against a \texttt{G21} quadrature rule with constant subinterval length in Figure~\ref{fig:quad_b}.
Both the time-stepping model with \texttt{K21} quadrature as well as the model with \texttt{G21} quadrature (which has a higher order than \texttt{K21}) are significantly more efficient, hence showing that the learned adaptive time-stepping is actually advantageous. The model with \texttt{G21} quadrature achieved an average error of about \num{6d-8} using 80.27 function evaluations per time while a constant time-stepper needs 87.20 function evaluations per time to achieve the same average error.
Thus, the adaptive time-stepping of our model leads to a reduction of required function evaluations of about 8\% compared to subintervals of constant size.

By optimizing the weights of the \texttt{G21} quadrature rule (as described in Section~\ref{subsec:optim_weights}) the performance can be enhanced further, leading to a reduction of required function evaluations of about 8.5\%.

\subsection{Numerical solution of differential equations}
\label{subsec:ResultsODE}

In the following, we apply our framework to several dynamical systems which exhibit chaotic motions. In this context, classical adaptive schemes often suffer from step rejections when the estimated local errors become too big. This drives up the number of function evaluations.

In all following experiments the integration scheme used by our time-stepping model is the same 5-th order Runge-Kutta scheme that is employed by the conventional \texttt{RK45} adaptive integrator, in order to enable a fair performance comparison.
For \texttt{RK45}, we use the implementation in the Python package SciPy.
We will show that our learned method obtains the required accuracy without the need of a step rejection mechanism and its performance exceeds \texttt{RK45}, even if function evaluations caused by rejected steps are not counted. This suggests that our approach learns to avoid the mistakes of traditional time-steppers.

\subsubsection{Example: Lorenz system}
    A frequently used benchmark for ODE integrators is the chaotic Lorenz system, i.\,e.,
    \begin{align*}
        \frac{d}{dt} \begin{pmatrix}
    				x_1\\
    				x_2\\
    				x_3
    			\end{pmatrix} = \begin{pmatrix}
    			\sigma (x_2 - x_1)\\
    			x_1 (\rho - x_3) - x_2\\
    			x_1 x_2 - \beta x_3
    		\end{pmatrix},
    \end{align*}
    with the standard parameter values $\sigma = 10$, $\beta = \frac{8}{3}$ and $\rho = 28$. In this regime, the system exhibits chaotic behavior.
    We train a model with step sizes 
    \[
        h \in \{0.02, 0.022, 0.025, 0.029, 0.033, 0.039, 0.045, 0.052, 0.060, 0.070\},
    \]
    a tolerance of $\tol=\num{d-4}$, and the time horizon $[0,T] = [0,100]$.
    For every training episode we sample a new initial condition $x(0)$ uniformly from the set $[-10, 10] \times [-10, 10] \times [15, 35]$, which contains a large portion of the attractor but also points away from the attractor.
    
\begin{figure}
    	\centering
    	\begin{subfigure}[b]{0.49\textwidth}
    	\includegraphics[width=\textwidth]{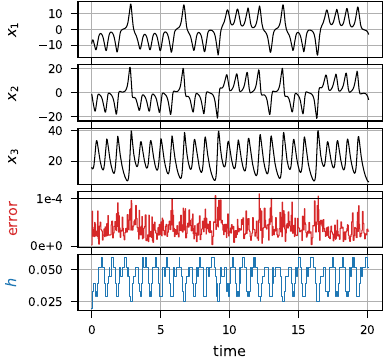}
    	\caption{time-stepping}
    	\label{fig:lorenz_a}
    	\end{subfigure}
    	\begin{subfigure}[b]{0.49\textwidth} \includegraphics[width=\textwidth]{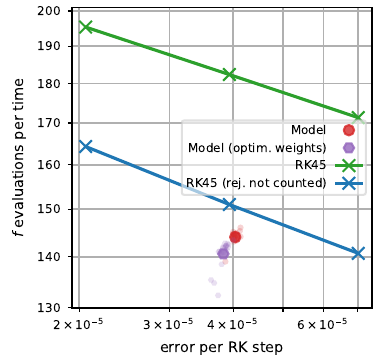}
    	\caption{Performance comparison}
    	\label{fig:lorenz_b}
    	\end{subfigure}
    	\caption{Lorenz System. (a) time-stepping obtained by the model. (b) Performance comparison between the model and \texttt{RK45}. For the time-stepping model (red) and the model with optimized weights (magenta) we plot the performance for an ensemble of 20 uniformly drawn initial conditions, integrated until $T=100$ (small dots), and the ensemble average (big dot). The data for \texttt{RK45} was obtained by choosing different desired tolerances for the local error estimates and averaging the performance over the ensemble. The blue data points (``rej. not counted'') show the performance if the function evaluations of step rejections of \texttt{RK45} are not counted.}
    	\label{fig:lorenz}
    \end{figure}

    The performance of the trained model is shown in Figure \ref{fig:lorenz}. We see in \ref{fig:lorenz_b} that it significantly outperforms \texttt{RK45}.
    Using on average of \num{143.9} function evaluations per time unit, the model achieves an average error of about \num{4d-5}, while \texttt{RK45} requires approximately \num{181.8} function evaluations to obtain the same average error.
    This corresponds to a reduction of the required function evaluations by approximately \SI{21}{\percent} at the same level of accuracy.
    As can be seen in Figure \ref{fig:lorenz_b}, the main reason for this substantial improvement is that \texttt{RK45} rejects poorly chosen time steps if the internal error estimation exceeds a certain bound. However, even without counting the rejections, the performance of our approach is slightly superior.
    By employing the technique of optimizing the weights of the used Runge-Kutta scheme (as described in \ref{subsec:optim_weights}), the performance can be enhanced further, so that the model can achieve the same accuracy as \texttt{RK45} with \SI{23}{\percent} fewer function evaluations.
    
\subsubsection{Example: Forced Van der Pol oscillator}
    \begin{figure}[t]
    	\centering
    	\begin{subfigure}[b]{0.49\textwidth}
    	\centering
    	\includegraphics[width=\linewidth]{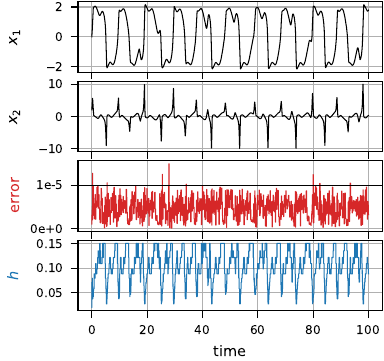}
    	\caption{time-stepping}
    	\label{fig:vdp_a}
    	\end{subfigure}
    	\hfill
    	\begin{subfigure}[b]{0.49\textwidth}
    	\centering
        \includegraphics[width=\linewidth]{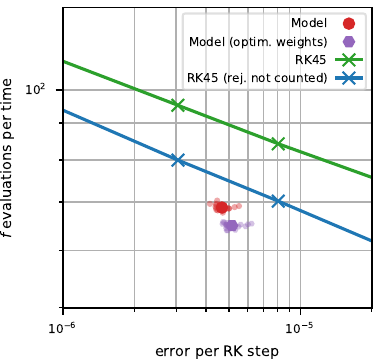}
    	\caption{Performance comparison}
    	\label{fig:vdp_b}
    	\end{subfigure}
    	\caption{Chaotic forced Van der Pol oscillator. (a) time-stepping obtained by the model. (b) Performance comparison between the model and \texttt{RK45}. For the time-stepping model (red) and the model with optimized weights (magenta) we plot the performance for an ensemble of 20 uniformly drawn initial conditions, integrated until $T=100$ (small dots), and the ensemble average (big dot). The data for \texttt{RK45} was obtained by choosing different desired tolerances for the local error estimates and averaging the performance over the ensemble. The blue data points (``rej. not counted'') show the performance if the function evaluations of step rejections of \texttt{RK45} are not counted.}
    	\label{fig:vdp}
    \end{figure}
    Next we study the Van der Pol oscillator with periodic forcing:
    \begin{align*}
        \frac{d^2 x}{dt^2} - \mu (1 - x^2) \frac{dx}{dt} + x = A \sin(\omega t)
    \end{align*}
    For the parameters $\mu = 5$, $A=5$ and $\omega=2.465$,
    proposed by Tsatsos \cite{tsatsos2008}, the system exhibits chaotic behavior.
    We train a time-stepper with $20$ step sizes evenly spaced on a $\log$ scale
    $h \in [0.02, 0.15]$ and $\tol = 10^{-5}$.
    
    Figure~\ref{fig:vdp_a} displays the time-stepping used by our model. Figure~\ref{fig:vdp_b} compares the efficiency of our model with SciPy's {\tt \texttt{RK45}}. While \texttt{RK45} suffers from step size rejections, which drives up the number of function evaluations, our model satisfies a given local error tolerance with fewer function evaluations.
    For an average error of about \num{5d-6}, it is using on average \num{68.8} function evaluations per time unit, while \texttt{RK45} needs approximately \num{90.7} function evaluations for the same average error.
    Our model needs approximately \SI{24}{\percent} less function evaluations than \texttt{RK45} with the same
    average error tolerance (\SI{28}{\percent} in case of the model with weight optimization enabled).
    Plotting the number of function evaluations of {\tt \texttt{RK45}} without counting function evaluations of rejected steps, we see a comparable behaviour to our method. The experiment suggests that our model learns to satisfy a local error tolerance without the need for step size rejections. This leads to lower computational costs compared to classical time-stepping control mechanisms.

    
    
    
    
    
    

\subsubsection{Example: Double Pendulum}
\label{subsec:example_doublependulum}
Another well-known chaotic system is the double pendulum, i.\,e., a pendulum with a second pendulum attached to its end.
The system is parametrized by the length and mass of each pendulum. (The limbs are modeled as massless and the masses of each pendulum as pointmasses at the end of the limbs.)
We set the two lengths and two masses to 1, gravitational acceleration to 10, and consider motion in the two-dimensional vertical plane only (pendulums can move left and right).
The dynamics of the double pendulum can then be characterized by an ODE of the state
\begin{align*}
    x = [\theta_1, \dot{\theta}_1, \theta_2, \dot{\theta}_2],
\end{align*}
where the $\theta_i$ refer to the angles of the pendulums and the $\dot{\theta}_i$ to the angular velocities \cite{shinbrot1992}.

We train a time-stepping model with 20 step sizes evenly spaced on a log scale $h \in [0.014, 0.1]$ and $\tol = 10^{-4}$.
The initial condition $x(0)$ is randomly drawn such that the resulting energy is always identical. We picked an energy level high enough for both pendulums to regularly flip.

\begin{figure}
    	\centering
    	\begin{subfigure}[b]{0.49\textwidth}
    	\includegraphics[width=\textwidth]{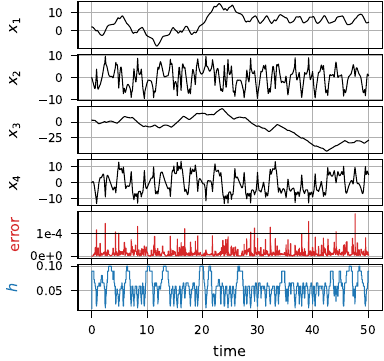}
    	\caption{time-stepping}
    	\label{fig:doublependulum_a}
    	\end{subfigure}
    	\begin{subfigure}[b]{0.49\textwidth} \includegraphics[width=\textwidth]{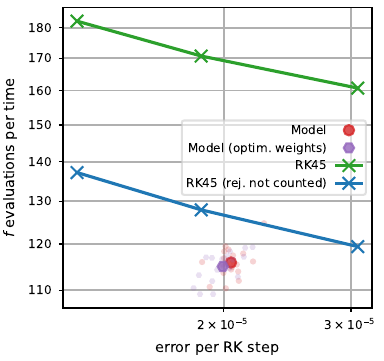}
    	\caption{Performance comparison}
    	\label{fig:doublependulum_b}
    	\end{subfigure}
    	\caption{Double Pendulum. (a) time-stepping obtained by the model, $x = [\theta_1, \dot{\theta}_1, \theta_2, \dot{\theta}_2]$ (b) Performance comparison between the model and \texttt{RK45}. For the time-stepping model (red) and the model with optimized weights (magenta) we plot the performance for an ensemble of 20 uniformly drawn initial conditions, integrated until $t_1=100$ (small dots), and the ensemble average (big dot). The data for \texttt{RK45} was obtained by choosing different desired tolerances for the local error estimates and averaging the performance over the ensemble. The blue data points (``rej. not counted'') show the performance if the function evaluations of step rejections of \texttt{RK45} are not counted.}
    	\label{fig:doublependulum}
\end{figure}

The trained model is shown in Figure \ref{fig:doublependulum}, and we see in \ref{fig:doublependulum_b} that it significantly outperforms \texttt{RK45}.
Using on average \num{115.9} function evaluations per time unit, the model achieves an average error of about \num{2d-5}, while \texttt{RK45} requires approximately \num{168.7} function evaluations to obtain the same average error.
This corresponds to a reduction of the required function evaluations by approximately \SI{31}{\percent} at the same level of accuracy.
As can be seen in Figure \ref{fig:doublependulum_b}, the main reason for this substantial improvement is that \texttt{RK45} rejects poorly chosen time steps if the internal error estimation exceeds a certain bound. However, even without counting the rejections, the performance of our approach is slightly superior.
By employing the technique of optimizing the weights of the used Runge-Kutta scheme (as described in \ref{subsec:optim_weights}), the performance can be enhanced slightly further, so that the model can achieve the same accuracy as \texttt{RK45} while using \SI{32}{\percent} less function evaluations.

\subsubsection{Example: H\'enon-Heiles System}
The H\'enon-Heiles system is related to galactic dynamics and exhibits chaotic behavior.
It is given by the ODE
\begin{align*}
        \frac{d}{dt} \begin{pmatrix}
    				x\\
    				p_x\\
    				y\\
    				p_y
    			\end{pmatrix} = \begin{pmatrix}
    			p_x \\
    			-x - 2xy\\
    			p_y \\
    			-y - (x^2 - y^2)
    		\end{pmatrix},
\end{align*}
where $(x, y)$ refer to coordinates in two-dimensional space and $(p_x, p_y)$ to the respective momenta.
We train a time-stepping model with step sizes
\begin{align*}
    h \in \{0.56, 0.58, 0.6, 0.62, 0.65, 0.68, 0.71, 0.74, 0.77, 0.8\}
\end{align*}
and $\tol = 10^{-4}$ by integrating from $t=0$ to $t=500$.
The initial positions $(x(0), y(0))$ are drawn uniformly from the triangle given by the contour-line of the H\'enon-Heiles potential with potential energy $\frac{1}{6}$ and the momenta such that the resulting total energy is equal to $\frac{1}{6}$.
As a consequence, the position $(x(t), y(t))$ is bounded within the above-mentioned triangle.

\begin{figure}
    	\centering
    	\begin{subfigure}[b]{0.49\textwidth}
    	\includegraphics[width=\textwidth]{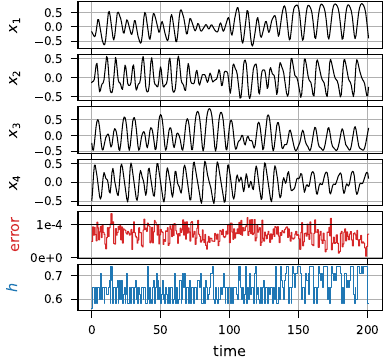}
    	\caption{time-stepping}
    	\label{fig:henon_a}
    	\end{subfigure}
    	\begin{subfigure}[b]{0.49\textwidth} \includegraphics[width=\textwidth]{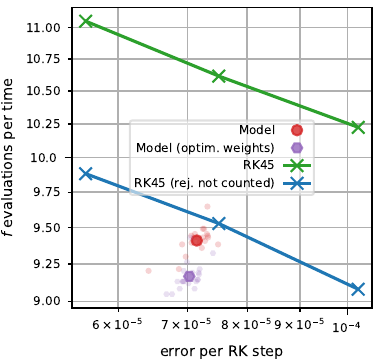}
    	\caption{Performance comparison}
    	\label{fig:henon_b}
    	\end{subfigure}
    	\caption{Henon-Heiles system. (a) time-stepping obtained by the model (b) Performance comparison between the model and \texttt{RK45}. For the time-stepping model (red) and the model with optimized weights (magenta) we plot the performance for an ensemble of 20 uniformly drawn initial conditions, integrated until $T=500$ (small dots), and the ensemble average (big dot). The data for \texttt{RK45} was obtained by choosing different desired tolerances for the local error estimates and averaging the performance over the ensemble. The blue data points (``rej. not counted'') show the performance if the function evaluations of step rejections of \texttt{RK45} are not counted.}
    	\label{fig:henon}
    \end{figure}

The trained model is shown in Figure \ref{fig:henon}, and we see in \ref{fig:henon_b} that it significantly outperforms \texttt{RK45}.
Using on average \num{9.41} function evaluations per time unit, the model achieves an average error of about \num{7d-5}, while \texttt{RK45} requires approximately \num{10.69} function evaluations to obtain the same average error.
This corresponds to a reduction of the required function evaluations by approximately \SI{7}{\percent} at the same level of accuracy.
As can be seen in Figure \ref{fig:henon_b}, the main reason for this substantial improvement is that \texttt{RK45} rejects poorly chosen time steps if the internal error estimation exceeds a certain bound. However, even without counting the rejections, the performance of our approach is slightly superior.
By employing the technique of optimizing the weights of the used Runge-Kutta scheme (as described in \ref{subsec:optim_weights}), the performance can be enhanced further, so that the model can achieve the same accuracy as \texttt{RK45} but with \SI{14}{\percent} fewer function evaluations.

\section{Conclusion} \label{sec:conclusion}
Combining classical numerical integrators with learned adaptive time-stepping strategies yields new numerical schemes which efficiently perform quadrature tasks and integrate differential equations. 
These schemes outperform state-of-the-art numerical procedures on problem classes for which traditional adaptive schemes show inefficient behaviour, such as chaotic systems. This is achieved by tailoring the step size controller to the problem classes at hand and to the desired tolerances.
Thus, our integration method combines the benefits of traditional numerical integration with data-driven optimization.
Our strategy is especially useful in situations, where a time-stepper can be trained offline and all computational tasks that need to be performed online belong to a restricted class of problems.



Experiments on the optimal choice of not only the time-stepping but also the associated quadrature or Runge-Kutta weights suggest a great potential for further improvement. Next to improvements in local accuracy, tailored weights can lead to structure preserving properties of the numerical scheme, which benefit long-term simulations.
For future work, it is of interest to develop a framework for a simultaneous selection of both time steps and quadrature weights and to design integrators that balance structure preservation, efficiency, and accuracy automatically.
Finally, the approach of a data-driven time-stepping strategy could be extended to a data-driven creation of space-time grids for the numerical integration of partial differential equations.

\section*{Source code} Our source code is freely accessible on GitHub: \url{https://github.com/lueckem/quadrature-ML}.

\section*{Author Contribution Statement}\mbox{}\\
Conceptualization and methodology: all authors.\\
Software: M. Lücke, K. Pfannenschmidt.\\
Writing: M. Lücke, C. Offen, S. Peitz, K. Pfannschmidt.

\bibliographystyle{abbrv}
\bibliography{main}
\end{document}